\documentclass[12pt,twoside]{amsart}
\usepackage{amssymb,amsmath,amsthm, amscd, enumerate, mathrsfs}
\usepackage{graphicx, hhline}
\usepackage[all]{xy}
\usepackage[usenames]{color}
\usepackage{hyperref}
\hypersetup{colorlinks=true}

\dedicatory{Dedicated to Professor Yujiro Kawamata on 
the occasion of his sixtieth birthday}

\title{Injectivity theorems}
\author{Osamu Fujino} 
\date{2015/7/3, version 1.60}
\address{Department of Mathematics, Graduate School of Science, 
Kyoto University, Kyoto 606-8502, Japan}
\email{fujino@math.kyoto-u.ac.jp}

\newcommand{\Spec}[0]{{\operatorname{Spec}}}

\newcommand{\Exc}[0]{{\operatorname{Exc}}}
\newcommand{\Supp}[0]{{\operatorname{Supp}}}
\newcommand{\Div}[0]{{\operatorname{Div}}}
\newcommand{\PerDiv}[0]{{\operatorname{PerDiv}}}
\newcommand{\Weil}[0]{{\operatorname{Weil}}}
\newcommand{\Gr}[0]{{\operatorname{Gr}}}
\newcommand{\Nlc}[0]{{\operatorname{Nlc}}}
\newtheorem{thm}{Theorem}[section]
\newtheorem{thm-sub}{Theorem}[subsection]
\newtheorem{lem}[thm]{Lemma}
\newtheorem{cor}[thm]{Corollary}
\newtheorem{prop}[thm]{Proposition}

\theoremstyle{definition}

\newtheorem{defn}[thm]{Definition}
\newtheorem{defn-sub}[thm-sub]{Definition}
\newtheorem{rem}[thm]{Remark}
\newtheorem{rem-sub}[thm-sub]{Remark}
\newtheorem*{ack}{Acknowledgments} 
\newtheorem{say}[thm]{}

\begin{document}

\maketitle 

\begin{abstract}
We prove some injectivity theorems. Our proof 
depends on the theory of mixed Hodge structures on cohomology 
with compact support. Our injectivity theorems would play crucial roles 
in the minimal model theory for 
higher-dimensional algebraic varieties. 
We also treat some applications. 
\end{abstract}

\tableofcontents 

\section{Introduction}\label{f-sec1} 

The following theorem is the main theorem of this paper, which is 
a slight generalization of \cite[Proposition 2.23]{book} (see also 
\cite{fujino-mi} and \cite[Theorem 3.1]{fujino-vani}) and 
is inspired by the main theorem of \cite{ambro}. We note that 
there are many contributors to this kind of injectivity theorem, for example, 
Tankeev, Koll\'ar, Esnault--Viehweg, Ambro, Fujino, and others. 

\begin{thm}[Main theorem]\label{f-thm1.1}
Let $X$ be a proper simple normal crossing algebraic variety and let $\Delta$ be 
an $\mathbb R$-Cartier $\mathbb R$-divisor 
on $X$ such that $\Supp \Delta$ is a simple normal crossing divisor on $X$ 
and that $\Delta$ is a boundary $\mathbb R$-divisor on $X$. 
Let $L$ be 
a Cartier divisor on $X$ and let $D$ be an effective Weil divisor on $X$ whose support 
is contained in $\Supp \Delta$. 
Assume that 
$L\sim _{\mathbb R}K_X+\Delta$. 
Then the natural homomorphism 
$$
H^q(X, \mathcal O_X(L))\to H^q(X, \mathcal O_X(L+D))
$$ 
induced by the inclusion $\mathcal O_X\to \mathcal O_X(D)$ is injective for every $q$. 
\end{thm}

\begin{rem}\label{f-rem1.2} 
In \cite[Proposition 2.23]{book}, the support of $D$ is 
assumed to be contained in $\Supp \{\Delta\}$, 
where $\{\Delta\}$ is the fractional part of $\Delta$. 
\end{rem}

\begin{rem}\label{f-rem1.3}
We will prove the relative version of Theorem \ref{f-thm1.1} in Theorem \ref{f-thm6.1}. 
The proof of Theorem \ref{f-thm6.1} uses \cite{bierstone-p}. 
Therefore, Theorem \ref{f-thm6.1} is a nontrivial generalization of Theorem \ref{f-thm1.1}. 
\end{rem}

We note that Theorem \ref{f-thm1.1} contains Theorem \ref{f-thm1.4}, which is 
equivalent to the main theorem of \cite{ambro} (see \cite[Theorem 2.3]{ambro}). 
Theorem \ref{f-thm1.4} shows that the notion of {\em{maximal non-lc ideal 
sheaves}} introduced in \cite{fst} 
is useful and has some nontrivial applications. 
For the details, see Section \ref{f-sec5}. 

\begin{thm}\label{f-thm1.4}
Let $X$ be a proper smooth algebraic variety 
and let $\Delta$ be a boundary $\mathbb R$-divisor 
on $X$ such that $\Supp \Delta$ is a simple normal 
crossing divisor on $X$. 
Let $L$ be a Cartier divisor on $X$ and let $D$ be an 
effective Cartier divisor on $X$ whose support 
is contained in $\Supp \Delta$. 
Assume that 
$L\sim _{\mathbb R}K_X+\Delta$. 
Then the natural homomorphism 
$$
H^q(X, \mathcal O_X(L))\to H^q(X, \mathcal O_X(L+D))
$$ 
induced by the inclusion $\mathcal O_X\to \mathcal O_X(D)$ is injective for every $q$. 
\end{thm}

A special case of Theorem \ref{f-thm1.1} implies a very powerful 
vanishing and torsion-free theorem for simple normal crossing 
pairs (see \cite[Theorem 1.1]{fujino-vani}). See also \cite{fujino-mi} and  
\cite[Theorem 2.38 and Theorem 2.39]{book}. 
It plays crucial roles for the study of {\em{semi-log canonical pairs}} and 
{\em{quasi-log varieties}} (see, \cite{book}, \cite{fujino-qlc}, \cite{fujino-slc}, 
and \cite{fuji-fuji}). 
 
More precisely, we obtain the following injectivity theorem for simple normal crossing 
pairs by using a special case of Theorem \ref{f-thm1.1}. 
 
\begin{thm}[{see \cite[Theorem 3.4]{fujino-vani}}]\label{f-thm1.5} 
Let $(X, \Delta)$ be a simple normal crossing 
pair such that $X$ is a proper algebraic variety and that $\Delta$ is a boundary 
$\mathbb R$-divisor on $X$. 
Let $L$ be a Cartier 
divisor on $X$ and let $D$ be an effective Cartier divisor that is permissible 
with 
respect to $(X, \Delta)$. 
Assume the following conditions:  
\begin{itemize}
\item[(i)] $L\sim _{\mathbb R}K_X+\Delta+H$, 
\item[(ii)] $H$ is a semi-ample $\mathbb R$-divisor, and 
\item[(iii)] $tH\sim _{\mathbb R} D+D'$ for some 
positive real number $t$, where 
$D'$ is an effective $\mathbb R$-Cartier 
$\mathbb R$-divisor that is permissible with respect to $(X, \Delta)$. 
\end{itemize}
Then the homomorphism 
$$
H^q(X, \mathcal O_X(L))\to H^q(X, \mathcal O_X(L+D)), 
$$ 
which is induced by the natural inclusion 
$\mathcal O_X\to \mathcal O_X(D)$, 
is injective for every $q$. 
\end{thm}

As an application of Theorem \ref{f-thm1.5}, we obtain Theorem \ref{f-thm1.6}, 
which is very important for the study of higher-dimensional algebraic varieties. 

\begin{thm}[{see \cite[Theorem 1.1]{fujino-vani}}]\label{f-thm1.6}
Let $(Y, \Delta)$ be a simple normal crossing 
pair such that $\Delta$ is a boundary $\mathbb R$-divisor on $Y$. 
Let $f:Y\to X$ be a proper morphism between algebraic varieties and let $L$ be a Cartier 
divisor on $Y$ such that 
$L-(K_Y+\Delta)$ is $f$-semi-ample. Let 
$q$ be an arbitrary non-negative integer. 
Then we have the following properties. 
\begin{itemize}
\item[(i)] Every associated prime of $R^qf_*\mathcal O_Y(L)$ is the 
generic point of the $f$-image of some stratum of $(Y, \Delta)$. 
\item[(ii)] 
Let $\pi:X\to V$ be a projective morphism to an algebraic variety $V$ such 
that $$L-(K_Y+\Delta)\sim _{\mathbb R}f^*H$$ for 
some $\pi$-ample $\mathbb R$-divisor $H$ on $X$. 
Then $R^qf_*\mathcal O_Y(L)$ is $\pi_*$-acyclic, that is, 
$$R^p\pi_*R^qf_*\mathcal O_Y(L)=0$$ for every $p>0$. 
\end{itemize}
\end{thm}
 
In this paper, we do not prove Theorem \ref{f-thm1.5} and Theorem \ref{f-thm1.6}.  
We only treat Theorem \ref{f-thm1.1} and Theorem \ref{f-thm1.4}. 
For the details of Theorem \ref{f-thm1.5} and Theorem \ref{f-thm1.6}, we 
recommend the reader to see \cite{fujino-vani}. 
 
Here, we quickly explain the main idea of the proof. 
 
\begin{say}[Idea of the proof]\label{f-say1.7} 
We give a proof of Theorem \ref{f-thm1.4} under the assumption that 
$\Delta$ is reduced and that $L\sim K_X+\Delta$.  

It is well-known that 
$$
E^{p,q}_1=H^q(X, \Omega^p_X(\log \Delta)\otimes \mathcal O_X(-\Delta))\Rightarrow 
H^{p+q}_c(X- \Delta, \mathbb C)
$$
degenerates at $E_1$ by Deligne's theory of mixed Hodge structures. 
This implies that the natural inclusion 
$$
\iota_!\mathbb C_{X- \Delta}\subset 
\mathcal O_X(-\Delta), 
$$ where $\iota:X- \Delta\to X$, induces surjections 
$$
\varphi^i:H^i(X, \iota_!\mathbb C_{X- \Delta})\to H^i(X, \mathcal O_X(-\Delta))
$$ 
for all $i$. On the other hand, 
we can easily see that 
$$
\iota_!\mathbb C_{X- \Delta}\subset \mathcal O_X(-\Delta-D)\subset \mathcal O_X(-\Delta)  
$$ because 
$\Supp D\subset \Supp \Delta$. 
Thus $\varphi^i$ factors as 
$$
H^i(X, \iota_!\mathbb C_{X- \Delta})
\to H^i(X, \mathcal O_X(-\Delta-D))
\to H^i(X, \mathcal O_X(-\Delta))
$$ 
for every $i$. Hence 
$$
H^i(X, \mathcal O_X(-\Delta-D))\to H^i(X, \mathcal O_X(-\Delta))
$$
is surjective for every $i$. By Serre duality, 
we obtain that 
$$
H^q(X, \mathcal O_X(K_X+\Delta))\to H^q(X, \mathcal O_X(K_X+\Delta+D))
$$ 
is injective for every $q$. 
\end{say}

In this paper, we use the notion of Du Bois complexes and Du Bois singularities for the 
proof of Theorem \ref{f-thm1.1} and Theorem \ref{f-thm1.4}. 
More precisely, we use the notion of Du Bois complexes for pairs, which 
is related to the mixed Hodge structures on cohomology with 
compact support. 
Consequently, the proof of Theorem \ref{f-thm1.1} is simpler 
than the arguments in \cite[Section 2.3 and Section 2.4]{book} 
(see also Section 3 and Section 4 in \cite{fujino-mi}). 
Note that we just need the $E_1$-degeneration of Hodge to de Rham type 
spectral sequences associated to 
the mixed Hodge structures on cohomology with compact support. 
We do not need the explicit descriptions of the weight filtrations. 

We strongly recommend the reader to see \cite{fujino-vani}. 
This paper and \cite{fujino-vani} simplify and 
generalize the main part of \cite[Chapter 2]{book} 
(see also \cite[Sections 3, 4, and 5]{fujino-mi}). 
We note that the foundation of the theory of semi-log canonical 
pairs discussed in \cite{fujino-slc} 
is composed of the results established in this paper 
and \cite{fujino-vani} (see \cite{fujino-mi} and \cite{book}). 

We summarize the contents of this paper. 
In Section \ref{f-sec2}, we collect some basic definitions and notations. 
In Section \ref{f-sec3}, we briefly review 
Du Bois complexes and Du Bois singularities. 
Section \ref{f-sec4} is devoted to the proof of Theorem \ref{f-thm1.1} and 
Theorem \ref{f-thm1.4}. In Section \ref{f-sec5}, we collect some miscellaneous comments 
on related topics, for example, Ambro's proof of 
the injectivity theorem in \cite{ambro}, 
the extension theorem from log canonical centers, etc. 
We also explain some interesting applications of Theorem \ref{f-thm1.4} 
due to Ambro (\cite{ambro}) in order to show how to use Theorem \ref{f-thm1.4}. 
In Section \ref{f-sec6}, we 
discuss the relative version of the main theorem:~Theorem 
\ref{f-thm6.1}. We also discuss some applications. 

\begin{ack} 
The author was partially supported by 
Grant-in-Aid for Young Scientists (A) 24684002 
from JSPS. 
He would like to thank Professors Akira Fujiki and 
Taro Fujisawa for answering his questions. 
He also would like to thank Professor Morihiko Saito.
The discussions with him on \cite{ffs} helped the author remove some ambiguities 
in a preliminary version of this paper. 
Finally, he thanks Professor Shunsuke Takagi for useful comments.  
\end{ack} 

We will work over $\mathbb C$, the field of complex numbers, throughout this paper. 
In this paper, a variety means a (not necessarily 
equidimensional) reduced separated scheme 
of finite type over $\mathbb C$. 
We will make use of the standard notation of the 
minimal model program as in \cite{fujino-fund}.  

\section{Preliminaries}\label{f-sec2}
First, we briefly recall basic definitions of 
divisors. We note that we have to 
deal with reducible varieties in this paper. 
For the details, see, for example, \cite[Section 2]{hartshorne} and \cite[Section 7.1]{liu}.  

\begin{say}\label{f-say2.1}
Let $X$ be a noetherian scheme with structure sheaf $\mathcal O_X$ 
and let $\mathcal K_X$ be the sheaf of total quotient rings of $\mathcal O_X$. 
Let $\mathcal K^*_X$ denote the (multiplicative) sheaf of invertible 
elements in 
$\mathcal K_X$, 
and $\mathcal O^*_X$ the sheaf of invertible elements in $\mathcal O_X$. 
We note that $\mathcal O_X\subset \mathcal K_X$ and 
$\mathcal O^*_X\subset \mathcal K^*_X$. 
\end{say}

\begin{say}[Cartier, $\mathbb Q$-Cartier, and 
$\mathbb R$-Cartier divisors]\label{f-say2.2}
A {\em{Cartier divisor}} $D$ on $X$ is a 
global section of $\mathcal K^*_X/\mathcal O^*_X$, 
that is, $D$ is an element of 
$H^0(X, \mathcal K^*_X/\mathcal O^*_X)$. 
A {\em{$\mathbb Q$-Cartier divisor}} (resp.~{\em{$\mathbb R$-Cartier 
divisor}}) is an element of 
$H^0(X, \mathcal K^*_X/\mathcal O^*_X)\otimes _{\mathbb Z}\mathbb Q$ 
(resp.~$H^0(X, \mathcal K^*_X/\mathcal O^*_X)\otimes _{\mathbb Z}\mathbb R$).  
\end{say}

\begin{say}[Linear, $\mathbb Q$-linear, and $\mathbb R$-linear equivalence]
\label{f-say2.3} 
Let $D_1$ and $D_2$ be two $\mathbb R$-Cartier divisors on $X$. 
Then $D_1$ is {\em{linearly}} (resp.~{\em{$\mathbb Q$-linearly}}, or 
{\em{$\mathbb R$-linearly}}) {\em{equivalent}} to $D_2$, denoted by 
$D_1\sim D_2$ (resp.~$D_1\sim_{\mathbb Q}D_2$, or 
$D_1\sim _{\mathbb R}D_2$) if $$D_1=D_2+\sum _{i=1}^k r_i(f_i)$$ such that 
$f_i\in \Gamma (X, \mathcal K^*_X)$ and $r_i\in \mathbb Z$ 
(resp.~$r_i\in \mathbb Q$, or $r_i\in \mathbb R$) for 
every $i$. We note that 
$(f_i)$ is a {\em{principal Cartier divisor}} 
associated to $f_i$, 
that is, the image of $f_i$ by $\Gamma (X, \mathcal K^*_X)\to 
\Gamma (X, \mathcal K^*_X/\mathcal O^*_X)$.  
\end{say}

\begin{say}[Supports]\label{f-say2.4}
Let $D$ be a Cartier divisor on $X$. 
The {\em{support}} of $D$, denoted by $\Supp D$, is 
the subset of $X$ consisting of points $x$ such that 
a local equation for $D$ is not in $\mathcal O^*_{X, x}$. The 
support of $D$ is a closed subset of $X$. 
\end{say}

\begin{say}[Weil divisors, $\mathbb Q$-divisors, and $\mathbb R$-divisors]
\label{f-say2.5} 
Let $X$ be an equidimensional reduced separated 
scheme of finite type over $\mathbb C$. We note that 
$X$ is not necessarily regular in codimension one. 
A {\em{\em{(}}Weil{\em{)}} divisor} 
$D$ on $X$ is a finite formal sum 
$$
\sum _{i=1}^nd_i D_i
$$
where $D_i$ is an irreducible reduced closed subscheme of $X$ of pure 
codimension one and $d_i$ is an integer for every $i$ such that 
$D_i\ne D_j$ for $i\ne j$.  

If $d_i\in \mathbb Q$ (resp.~$d_i\in \mathbb R$) for every $i$, 
then $D$ is called a {\em{$\mathbb Q$-divisor}} (resp.~{\em{$\mathbb R$-divisor}}). 
We define the {\em{round-up}} $\lceil D\rceil=\sum _{i=1}^r\lceil d_i\rceil D_i$ 
(resp.~the {\em{round-down}} $\lfloor D\rfloor =\sum_{i=1}^r\lfloor d_i\rfloor D_i$), 
where for every real number $x$, $\lceil x\rceil$ (resp.~$\lfloor x\rfloor$) is the integer 
defined by $x\leq \lceil x\rceil<x+1$ (resp.~$x-1<\lfloor x\rfloor \leq x$). 
The {\em{fractional part}} $\{D\}$ of 
$D$ denotes $D-\lfloor D\rfloor$. 
We call $D$ a {\em{boundary}} 
$\mathbb R$-divisor if $0\leq d_i\leq 1$ for every $i$. 

We put $$
D^{\leq k}=\sum _{d_i\leq k}d_i D_i, \quad D^{\geq k}=\sum _{d_i\geq k}d_i D_i, \quad 
D^{=k}=\sum _{d_i=k}d_i D_i=k\sum _{d_i=k}D_i
$$ 
and 
$$
{}^k\!D=\sum _{d_i=k}D_i 
$$
for every $k\in \mathbb R$. 
We note that $D^{=1}={}^1\!D$. 
\end{say}

Next, we recall the definition of {\em{simple normal crossing 
pairs}}. 

\begin{defn}[Simple normal crossing pairs]\label{f-def2.6} 
We say that the pair $(X, D)$ is {\em{simple normal crossing}} at 
a point $a\in X$ if $X$ has a Zariski open 
neighborhood $U$ of $a$ that can be embedded in a smooth 
variety 
$Y$, 
where $Y$ has regular system of parameters $(x_1, \cdots, x_p, y_1, \cdots, y_r)$ at 
$a=0$ in which $U$ is defined by a monomial equation 
$$
x_1\cdots x_p=0
$$ 
and $$
D=\sum _{i=1}^r \alpha_i(y_i=0)|_U, \quad  \alpha_i\in \mathbb R. 
$$ 
We say that $(X, D)$ is a {\em{simple 
normal crossing pair}} if it is simple normal crossing at every point of $X$. 
If $(X, 0)$ is a simple normal crossing pair, 
then $X$ is called a {\em{simple normal crossing 
variety}}. If $X$ is a simple normal crossing variety, 
then $X$ has only Gorenstein singularities. 
Thus, it has an invertible dualizing sheaf $\omega_X$. 
Therefore, we can define the {\em{canonical divisor $K_X$}} such that 
$\omega_X\simeq \mathcal O_X(K_X)$ (cf.~\cite[Section 7.1 Corollary 1.19]{liu}). 
It is a Cartier divisor on $X$ and is well-defined up to linear equivalence. 
\end{defn}

We note that a simple normal crossing pair is called a {\em{semi-snc pair}} 
in \cite[Definition 1.10]{kollar-book}. 

\begin{defn}[Strata and permissibility]\label{f-def2.7} 
Let $X$ be 
a simple normal crossing variety and let $X=\bigcup _{i\in I}X_i$ be the 
irreducible decomposition of $X$. 
A {\em{stratum}} of $X$ is an irreducible component of $X_{i_1}\cap 
\cdots \cap X_{i_k}$ for some 
$\{i_1, \cdots, i_k\}\subset I$. 
A Cartier divisor $D$ on $X$ is {\em{permissible}} if $D$ contains no 
strata of $X$ in its support. 
A finite $\mathbb Q$-linear (resp.~$\mathbb R$-linear) combination of 
permissible 
Cartier divisors is called a {\em{permissible 
$\mathbb Q$-divisor}} (resp.~{\em{$\mathbb R$-divisor}}) on 
$X$. 
\end{defn}

\begin{say}\label{f-say2.8}
Let $X$ be a simple normal crossing variety. 
Let $\PerDiv (X)$ be the abelian group generated by permissible 
Cartier divisors on $X$ and let $\Weil (X)$ be the abelian group 
generated by Weil divisors on $X$. 
Then we can define natural injective homomorphisms of abelian groups 
$$
\psi:\PerDiv (X)\otimes _{\mathbb Z}\mathbb K\to \Weil (X)\otimes _{\mathbb Z}\mathbb K 
$$ 
for $\mathbb K=\mathbb Z$, $\mathbb Q$, and $\mathbb R$. 
Let $\nu:\widetilde X\to X$ be the normalization. 
Then we have the following commutative diagram. 
$$
\xymatrix{
\Div (\widetilde X)\otimes _{\mathbb Z}\mathbb K 
\ar[r]^{\sim}_{\widetilde \psi}& \Weil (\widetilde X)\otimes 
_{\mathbb Z}\mathbb K\ar[d]^{\nu_*} \\
\PerDiv (X)\otimes _{\mathbb Z}\mathbb K\ar[r]_{\psi}
\ar[u]^{\nu^*}&\Weil (X)\otimes _{\mathbb Z}\mathbb K
}
$$ 
Note that $\Div (\widetilde X)$ is the abelian group generated by 
Cartier divisors on $\widetilde X$ and 
that $\widetilde \psi$ is an isomorphism 
since $\widetilde X$ is smooth. 

By $\psi$, every permissible divisor (resp.~$\mathbb Q$-divisor  
or $\mathbb R$-divisor) can be considered 
as a Weil divisor (resp.~$\mathbb Q$-divisor or $\mathbb R$-divisor). 
Therefore, various operations, 
for example, $\lfloor D\rfloor$, $\{D\}$, and so on, 
make sense for a permissible $\mathbb R$-divisor $D$ on $X$. 
\end{say}

\begin{defn}[Simple normal crossing divisors]\label{f-def2.9}
Let $X$ be a simple normal crossing variety and 
let $D$ be a Cartier divisor on $X$. 
If $(X, D)$ is a simple normal crossing pair and $D$ is reduced, 
then $D$ is called a {\em{simple normal crossing divisor}} on $X$. 
\end{defn}

\begin{rem}\label{f-rem2.10}
Let $X$ be a simple normal crossing variety and let $D$ be a 
$\mathbb K$-divisor on $X$ where 
$\mathbb K=\mathbb Q$ or $\mathbb R$. 
If $\Supp D$ is a simple normal crossing divisor on $X$ and $D$ is $\mathbb K$-Cartier, 
then $\lfloor D\rfloor$ and $\lceil D\rceil$ (resp.~$\{D\}$, $D^{<1}$, 
and so on) are Cartier (resp.~$\mathbb K$-Cartier) divisors 
on $X$ (cf.~\cite[Section 8]{bierstone-p}). 
\end{rem}

The following lemma is easy but important. 

\begin{lem}\label{f-lem2.11}
Let $X$ be a simple normal crossing variety and 
let $B$ be a permissible $\mathbb R$-divisor on 
$X$ such that $\lfloor B\rfloor=0$. 
Let $A$ be a Cartier divisor on $X$. 
Assume that 
$A\sim _{\mathbb R}B$. Then 
there exists a permissible $\mathbb Q$-divisor 
$C$ on $X$ such that $A\sim _{\mathbb Q}C$, $\lfloor C\rfloor =0$, and 
$\Supp C=\Supp B$. 
\end{lem}
\begin{proof}
We can write $B=A+\sum _{i=1}^k r_i (f_i)$, 
where $f_i\in \Gamma (X, \mathcal K^* _X)$ and 
$r_i\in \mathbb R$ for every $i$. 
Let $P\in X$ be a scheme theoretic point corresponding to 
some stratum of $X$. 
We consider the following affine map 
$$
\mathbb K^k\to H^0(X_P, 
\mathcal K^*_{X_P}/\mathcal O^{*}_{X_P})\otimes _{\mathbb Z} \mathbb K
$$ 
induced by $(a_1, \cdots, a_k)\mapsto A+\sum _{i=1}^k a_i (f_i)$, where 
$X_P=\Spec \mathcal O_{X, P}$ and $\mathbb K=\mathbb Q$ or $\mathbb R$. 
Then we can check that 
$$
\mathcal P=\{(a_1, \cdots, a_k)\in \mathbb R^k\,|\, A+\sum _i a_i (f_i)\ 
\text{is permissible}\}\subset \mathbb R^k
$$ 
is an affine subspace of $\mathbb R^k$ defined over $\mathbb Q$. 
Therefore, we see that 
$$
\mathcal S=\{(a_1, \cdots, a_k)\in 
\mathcal P\, |\, \Supp (A+\sum _i a_i (f_i))\subset \Supp B\}\subset 
\mathcal P 
$$ is an affine subspace of $\mathbb R^k$ defined over $\mathbb Q$. 
Since $(r_1, \cdots, r_k)\in \mathcal S$, we know 
that $\mathcal S\ne \emptyset$. 
We take a point $(s_1, \cdots, s_k)\in \mathcal S \cap \mathbb Q^k$ which 
is general in $\mathcal S$ and sufficiently close to 
$(r_1, \cdots, r_k)$ and put $C=A+\sum _{i=1}^k s_i (f_i)$. By construction, 
$C$ is a permissible $\mathbb Q$-divisor such that 
$C\sim _{\mathbb Q}A$, $\lfloor C\rfloor =0$, and 
$\Supp C=\Supp B$. 
\end{proof}

\section{A quick review of Du Bois complexes}\label{f-sec3}

In this section, we briefly review Du Bois complexes and Du Bois singularities. 
For the details, see, for example, \cite{dubois}, 
\cite{steenbrink}, \cite[Expos\'e V]{gnpp}, \cite{saito}, 
\cite{ps}, \cite{kovacs}, and \cite[Chapter 6]{kollar-book}.  

\begin{say}[Du Bois complexes]\label{f-say3.1} 
Let $X$  be an algebraic variety. 
Then we can associate a filtered complex $(\underline {\Omega}^\bullet_X, F)$ 
called the {\em{Du Bois complex}} of 
$X$ in a suitable derived category $D^b_{\mathrm{diff, coh}}(X)$ 
(see \cite[1.~Complexes filtr\'es d'op\'erateurs diff\'erentiels d'ordre $\leq 1$]{dubois}).  
We put 
$$
\underline {\Omega}^0_X=\Gr^0_F\underline {\Omega}^\bullet_X. 
$$ 
There is a natural map $(\Omega^\bullet_X, 
\sigma)\to (\underline {\Omega}^\bullet_X, F)$. 
It induces $\mathcal O_X\to \underline {\Omega}^0_X$. 
If $\mathcal O_X\to \underline {\Omega}^0_X$ is a quasi-isomorphism, 
then $X$ is said to have {\em{Du Bois singularities}}. 
We sometimes simply say that $X$ is {\em{Du Bois}}. 
Let $\Sigma$ be a reduced closed subvariety of $X$. Then there is a natural 
map $\rho:(\underline {\Omega}^\bullet_X, F)\to (\underline{\Omega}^\bullet_\Sigma, F)$ 
in $D^b_{\mathrm{diff, coh}}(X)$. 
By taking the cone of $\rho$ with a shift by one, 
we obtain a filtered 
complex $(\underline{\Omega}^\bullet_{X, \Sigma}, F)$ in $D^b_{{\mathrm{diff, coh}}}(X)$. 
Note that $(\underline {\Omega}^\bullet_{X, \Sigma}, F)$ was essentially introduced 
by Steenbrink in \cite[Section 3]{steenbrink}. 
We put 
$$
\underline {\Omega}^0_{X, \Sigma}=\Gr ^0_F\underline {\Omega}^\bullet_{X, \Sigma}. 
$$
Then there are a map $\mathcal J_\Sigma\to \underline {\Omega}^0_{X, \Sigma}$, 
where $\mathcal J_\Sigma$ is the defining ideal sheaf of $\Sigma$ on $X$, and 
the following commutative diagram 
\begin{equation*}
\xymatrix{\mathcal J_\Sigma\ar[r] \ar[d]& 
\mathcal O_X\ar[r]\ar[d]&\mathcal O_\Sigma\ar[r]^{+1}\ar[d]&\\
 \underline \Omega^0_{X, \Sigma} 
 \ar[r] & \underline \Omega^0_X\ar[r]& \underline {\Omega}^0_\Sigma\ar[r]^{+1}&}
\end{equation*}
in the derived category $D^b_{\mathrm{coh}}(X)$ 
(see also Remark \ref{f-rem3.3} below). 
\end{say}

By using the theory of mixed 
Hodge structures on cohomology with compact support, we 
have the following theorem. 

\begin{thm}\label{f-thm3.2} 
Let $X$ be a variety and 
let $\Sigma$ be a reduced closed subvariety of $X$. 
We put $j:X-\Sigma \hookrightarrow X$. 
Then we have the following 
properties. 
\begin{itemize}
\item[(1)] 
The complex $(\underline {\Omega}^{\bullet}_{X, \Sigma})^{\mathrm{an}}$ 
is a resolution of $j_!\mathbb C_{X^{\mathrm{an}}-\Sigma^{\mathrm{an}}}$. 
\item[(2)] If in addition $X$ is proper, then the spectral sequence 
$$
E^{p, q}_1=\mathbb H^q(X, \underline 
{\Omega}^{p}_{X, \Sigma})\Rightarrow H^{p+q}(X^{\mathrm{an}}, j_!\mathbb C_
{X^{\mathrm{an}}-\Sigma^{\mathrm{an}}})
$$
degenerates at $E_1$, where $\underline 
{\Omega}^p_{X, \Sigma}=\Gr^p_F \underline {\Omega}^\bullet_{X, \Sigma}
[p]$. 
\end{itemize}
\end{thm}

From now on, we will simply write $X$ (resp.~$\mathcal O_X$ and so on) 
to express $X^{\mathrm{an}}$ (resp.~$\mathcal O_{X^{\mathrm{an}}}$ and so on) 
if there is no risk of confusion. 

\begin{proof}Here, we use the formulation of \cite[\S3.3 and \S3.4]{ps}. 
We assume that $X$ is proper. 
We take cubical hyperresolutions $\pi_X:X_{\bullet}\to X$ and 
$\pi_\Sigma: \Sigma_{\bullet}\to \Sigma$ fitting in a commutative 
diagram. 
$$
\xymatrix{
\Sigma_{\bullet}\ar[d]_{\pi_\Sigma}\ar[r]&X_{\bullet}\ar[d]^{\pi_X}\\
\Sigma\ar[r]_{\iota}&X
}
$$
Let $\mathcal Hdg(X):=R\pi_{X*}\mathcal Hdg^\bullet(X_\bullet)$ 
be a mixed Hodge complex of sheaves on $X$ giving 
the natural mixed Hodge structure 
on $H^\bullet(X, \mathbb Z)$ (see \cite[Definition 5.32 and Theorem 5.33]{ps}). 
We can obtain a mixed Hodge complex 
of sheaves $\mathcal Hdg(\Sigma):=
R\pi_{\Sigma*}\mathcal Hdg^\bullet(\Sigma_\bullet)$ on $\Sigma$ 
analogously. Roughly speaking, by 
forgetting the weight filtration and the $\mathbb Q$-structure of 
$\mathcal Hdg(X)$ and considering it in $D^b_{\mathrm{diff, coh}}(X)$, 
we obtain the Du Bois complex $(\underline {\Omega}_X^\bullet, F)$ of $X$ 
(see \cite[Expos\'e V (3.3) Th\'eor\'eme]{gnpp}). 
We can also obtain the Du Bois complex $(\underline{\Omega}_\Sigma^\bullet, F)$ 
of $\Sigma$ analogously. 
By taking the mixed cone of $\mathcal Hdg(X)\to \iota_*\mathcal Hdg(\Sigma)$ with 
a shift by one, we obtain a mixed Hodge complex of sheaves on $X$ giving the natural 
mixed Hodge structure on $H_c^\bullet(X-\Sigma, \mathbb Z)$ 
(see \cite[5.5~Relative Cohomology]{ps}). 
Roughly speaking, by forgetting the weight filtration and the $\mathbb Q$-structure, 
we obtain the desired filtered complex $(\underline{\Omega}_{X, \Sigma}^\bullet, F)$ in 
$D^b_{\mathrm{diff, coh}}(X)$. 
When $X$ is not proper, we take completions of 
$\overline X$ and $\overline {\Sigma}$ of $X$ and $\Sigma$ and 
apply the above arguments to $\overline X$ and $\overline \Sigma$. 
Then we restrict everything to $X$. 
The properties (1) and (2) obviously hold by the above 
description of $(\underline{\Omega}_{X, \Sigma}^\bullet, F)$. 
By the above construction and description of $(\underline{\Omega}_{X, \Sigma}^\bullet, 
F)$, we know that the map $\mathcal J_{\Sigma}\to \underline{\Omega}_{X, \Sigma}^0$ 
in $D^b_{\mathrm{coh}}(X)$ is induced by natural 
maps of complexes. 
\end{proof}

\begin{rem}\label{f-rem3.3}
Note that the Du Bois complex $\underline{\Omega}_X^\bullet$ is nothing but 
the filtered complex $R\pi_{X*}(\Omega_{X_\bullet}^\bullet, F)$. 
For the details, see \cite[Expos\'e V (3.3) Th\'eor\'eme and 
(3.5) D\'efinition]{gnpp}. 
Therefore, the Du Bois complex of the pair $(X, \Sigma)$ 
is given by 
$$
\mathrm{Cone}^\bullet(R\pi_{X*}(\Omega_{X_\bullet}^\bullet, F)\to 
\iota_*R\pi_{\Sigma*}(\Omega_{\Sigma_\bullet}^\bullet, F))[-1]. 
$$ 
By the construction of $\underline {\Omega}_X^\bullet$, there is a natural map 
$a_X: \mathcal O_X\to \underline{\Omega}_X^\bullet$ which 
induces $\mathcal O_X\to \underline{\Omega}_X^0$ in $D_{\mathrm{coh}}^b(X)$. 
Moreover, the composition of $a_X^{\mathrm{an}}: 
\mathcal O_{X^{\mathrm{an}}}\to (\underline{\Omega}_X^\bullet)^{\mathrm{an}}$ 
with the natural inclusion $\mathbb C_{X^\mathrm{an}}\subset 
\mathcal O_{X^\mathrm{an}}$ induces a quasi-isomorphism 
$\mathbb C_{X^{\mathrm{an}}}\overset{\simeq}{\longrightarrow}
(\underline{\Omega}_X^\bullet)^{\mathrm{an}}$. 
We have a natural map 
$a_\Sigma:\mathcal O_\Sigma\to \underline {\Omega}_\Sigma^\bullet$ 
with the same properties as $a_X$ and 
the following commutative diagram. 
$$
\xymatrix{\mathcal O_X\ar[r] \ar[d]_{a_X}&\mathcal O_\Sigma\ar[d]^{a_\Sigma}\\
\underline {\Omega}_X^\bullet \ar[r]&\underline {\Omega}_\Sigma^\bullet
}
$$
Therefore, we have a natural map 
$b: \mathcal J_\Sigma\to \underline{\Omega}_{X, \Sigma}^\bullet$ such that 
$b$ induces 
$\mathcal J_\Sigma\to \underline{\Omega}_{X, \Sigma}^0$ in 
$D_{\mathrm{coh}}^b(X)$ and 
that the composition of $b^{\mathrm{an}}: (\mathcal J_\Sigma)^{\mathrm{an}}
\to (\underline{\Omega}_{X, \Sigma}^\bullet)^{\mathrm{an}}$ with 
the natural inclusion $j_!\mathbb C_{X^{\mathrm{an}}-\Sigma^{\mathrm{an}}}\subset 
(\mathcal J_\Sigma)^{\mathrm{an}}$ induces a quasi-isomorphism 
$j_!\mathbb C_{X^{\mathrm{an}}-\Sigma^{\mathrm{an}}}\overset{\simeq}{\longrightarrow} 
(\underline{\Omega}_{X, \Sigma}^\bullet)^{\mathrm{an}}$.  
We need the weight filtration and the $\mathbb Q$-structure 
in order to prove the $E_1$-degeneration of Hodge to 
de Rham type spectral sequence. 
We used the framework of \cite[\S3.3 and \S3.4]{ps} 
because we had to check that various diagrams related to 
comparison morphisms are commutative (see \cite[Remark 3.23]{ps}) 
for 
the proof of Theorem \ref{f-thm3.2} (2) and so on. 
\end{rem}

Let us recall the definition of {\em{Du Bois pairs}} by \cite[Definition 3.13]{kovacs}. 

\begin{defn}[Du Bois pairs]\label{f-def3.4} 
With the notation of \ref{f-say3.1} and Theorem \ref{f-thm3.2}, 
if the map $\mathcal J_\Sigma\to  \underline{\Omega}^0_{X, \Sigma}$ is 
a quasi-isomorphism, 
then the pair $(X, \Sigma)$ is called a {\em{Du Bois pair}}. 
\end{defn}

By the definitions, we can easily check the following useful proposition. 

\begin{prop}\label{f-prop3.5} 
With the notation of \ref{f-say3.1} and Theorem \ref{f-thm3.2}, 
we assume that both $X$ and $\Sigma$ are Du Bois. 
Then the pair $(X, \Sigma)$ is a 
Du Bois pair, that is, $\mathcal J_\Sigma\to \underline{\Omega}
^0_{X, \Sigma}$ is a quasi-isomorphism. 
\end{prop}

Let us recall the following well-known results on Du Bois singularities. 

\begin{thm}\label{f-thm3.6} 
Let $X$ be a normal algebraic variety with only quotient singularities. 
Then $X$ has only rational singularities. In particular, $X$ is Du Bois. 
\end{thm}

Theorem \ref{f-thm3.6} follows from, 
for example, \cite[5.2.~Th\'eor\`eme]{dubois}, and 
\cite{kovacs1}. 
Lemma \ref{f-lem3.7} will play an important role 
in the proof of Theorem \ref{f-thm1.4}. 

\begin{lem}\label{f-lem3.7} 
Let $X$ be a variety with 
closed subvarieties $X_1$ and $X_2$ such that $X=X_1\cup X_2$. 
Assume that $X_1$, $X_2$, and $X_1\cap X_2$ are 
Du Bois. Note that, in particular, we assume that 
$X_1\cap X_2$ is reduced. Then $X$ is Du Bois. 
\end{lem} 

For the proof of Lemma \ref{f-lem3.7}, see, for example, \cite[Lemma 3.4]{schwede}. 
We close this section with a remark on Du Bois singularities. 

\begin{rem}[Du Bois singularities and log canonical singularities]
Koll\'ar and Kov\'acs established that 
log canonical singularities are Du Bois in \cite{kk}. 
Moreover, semi-log canonical singularities are 
Du Bois (see \cite[Corollary 6.32]{kollar-book}). 
We note that the arguments in \cite{kk} heavily depend 
on the recent developments of the minimal model 
program by Birkar--Cascini--Hacon--M\textsuperscript{c}Kernan and 
the results by Ambro and Fujino (see, 
for example, \cite{ambro-qlog}, \cite{book}, \cite{fujino-non}, 
and \cite{fujino-fund}). 
We need a special case of Theorem \ref{f-thm1.6} for the arguments 
in \cite{kk}. 
In this paper, we will just use Du Bois 
complexes for cyclic covers of simple normal crossing 
pairs. 
Our proof in Section \ref{f-sec4} is independent of the deep result in \cite{kk}. 

The fact that (semi-)log canonical singularities are Du Bois does not seem to be so useful 
when we consider various Kodaira-type vanishing theorems for 
(semi-)log canonical 
pairs. 
This is because (semi-)log canonical 
singularities are not necessarily Cohen--Macaulay. 
The approach to various Kodaira-type vanishing theorems 
for semi-log canonical pairs in \cite{fujino-slc} 
is based on the vanishing theorem in \cite{fujino-vani} 
(see Theorem \ref{f-thm1.6}, \cite{fujino-mi}, and \cite{book}) 
and the theory of partial 
resolution of singularities for reducible varieties (see \cite{bierstone-p}). 
\end{rem}

\section{Proof of theorems}\label{f-sec4}

In this section, we prove Theorem \ref{f-thm1.1} and Theorem \ref{f-thm1.4}. 

\begin{proof}[Proof of Theorem \ref{f-thm1.4}]
Without loss of generality, we may assume 
that $X$ is connected. We set $S=\lfloor \Delta\rfloor$ and $B=\{\Delta\}$. 
By perturbing $B$, 
we may assume that $B$ is a $\mathbb Q$-divisor (cf.~Lemma \ref{f-lem2.11}).  
We set $\mathcal M=\mathcal O_X(L-K_X-S)$. 
Let $N$ be the smallest positive integer 
such that $N L\sim N(K_X+ S+ B)$. 
In particular, $N B$ is an integral Weil divisor. 
We take the $N$-fold cyclic cover 
$$\pi': Y'=\Spec_X\!\bigoplus _{i=0}^{N-1} \mathcal M^{-i}\to 
X$$ associated to 
the section $N B\in |\mathcal M^{N}|$. 
More precisely, let $s\in H^0(X, \mathcal M^{N})$ be a section 
whose zero divisor is $N B$. 
Then the dual of $s:\mathcal O_X\to \mathcal M^{N}$ 
defines an $\mathcal O_X$-algebra structure on 
$\bigoplus ^{N-1}_{i=0} \mathcal M^{-i}$. 
Let $Y\to Y'$ be the normalization and 
let $\pi:Y\to X$ be the composition morphism. 
It is well-known that 
$$
Y=\Spec _X\bigoplus _{i=0}^{N-1}\mathcal M^{-i}(\lfloor iB\rfloor). $$  
For the details, see \cite[3.5.~Cyclic covers]{ev}. 
Note that $Y$ has only quotient singularities. 
We set $T=\pi^*S$. 
Let $T=\sum _{i\in I}T_i$ be the irreducible decomposition. 
Then every irreducible component of $T_{i_1}\cap \cdots \cap T_{i_k}$ has only 
quotient singularities for every $\{i_1, \cdots, i_k\}\subset I$. 
Hence it is easy to see that both $Y$ and $T$ have only Du Bois singularities by 
Theorem \ref{f-thm3.6} and Lemma \ref{f-lem3.7} 
(see also \cite{ishida}). 
Therefore, the pair $(Y, T)$ is a Du Bois pair by Proposition \ref{f-prop3.5}. 
This means that $\mathcal O_Y(-T)\to 
\underline {\Omega}^0_{Y, T}$ is a quasi-isomorphism. 
See also \cite[3.4]{ffs}. We note that 
$T$ is Cartier. 
Hence $\mathcal O_Y(-T)$ is the defining ideal sheaf of $T$ on $Y$. 
The $E_1$-degeneration of 
$$
E^{p, q}_1=\mathbb 
H^q(Y, \underline {\Omega}^p_{Y, T})\Rightarrow H^{p+q}(Y, j_!\mathbb C_{Y-T})
$$ 
implies that the homomorphism 
$$
H^q(Y, j_!\mathbb C_{Y-T})\to H^q(Y, \mathcal O_Y(-T))
$$
induced by the natural inclusion 
$$
j_!\mathbb C_{Y-T}\subset \mathcal O_Y(-T)
$$ 
is surjective for every $q$ (see Remark \ref{f-rem3.3}). 
By taking a suitable direct summand $$\mathcal C\subset \mathcal M^{-1}(-S)$$ of 
$$\pi_*(j_!\mathbb C_{Y-T})\subset \pi_*\mathcal O_Y(-T), $$ we 
obtain a surjection 
$$
H^q(X, \mathcal C)\to H^q(X, \mathcal M^{-1}(-S)) 
$$ 
induced by the natural inclusion $\mathcal C\subset \mathcal M^{-1}(-S)$ for 
every $q$. We can check 
the following simple property by examining the monodromy 
action of the Galois group $\mathbb Z/N\mathbb Z$ of 
$\pi:Y\to X$ 
on $\mathcal C$ around $\Supp B$. 

\begin{lem}[{cf.~\cite[Corollary 2.54]{km}}]\label{f-lem4.1}
Let $U\subset X$ be 
a connected open set such that 
$U\cap \Supp \Delta\ne \emptyset$. 
Then $H^0(U, \mathcal C|_{U})=0$. 
\end{lem}
\begin{proof}
If $U\cap \Supp B\ne \emptyset$, then $H^0(U, \mathcal C|_U)=0$ since the monodromy 
action on $\mathcal C|_{U\setminus \Supp B}$ around $\Supp B$ is nontrivial. 
If $U\cap \Supp S\ne \emptyset$, then $H^0(U, \mathcal C|_U)=0$ since 
$\mathcal C$ is a direct summand of $\pi_*(j_!\mathbb C _{Y-T})$ and $T=\pi^*S$. 
\end{proof}

This property is utilized via the following fact. 
The proof is obvious. 

\begin{lem}[{cf.~\cite[Lemma 2.55]{km}}]\label{f-lem4.2}
Let $F$ be a sheaf of Abelian groups on a topological 
space 
$X$ and $F_1, F_2\subset F$ subsheaves. 
Let $Z\subset X$ be a closed subset. Assume 
that 
\begin{itemize}
\item[(1)] $F_2|_{X- Z}=F|_{X- Z}$, and 
\item[(2)] if $U$ is connected, open 
and $U\cap Z\ne \emptyset$, then 
$H^0(U, F_1|U)=0$. 
\end{itemize}
Then $F_1$ is a subsheaf of $F_2$. 
\end{lem}

As a corollary, we obtain: 

\begin{cor}[{cf.~\cite[Corollary 2.56]{km}}]\label{f-cor4.3}
Let $M\subset \mathcal M^{-1}(-S)$ be a subsheaf such that 
$M|_{X- \Supp \Delta}
=\mathcal M^{-1}(-S)|_{X- \Supp \Delta}$. 
Then the injection 
$$
\mathcal C\to \mathcal M^{-1}(-S) 
$$ 
factors as 
$$
\mathcal C \to M\to \mathcal M^{-1}(-S). 
$$ Therefore, 
$$H^q(X, M)\to H^q(X, \mathcal M^{-1}(-S))
$$ 
is surjective for every $q$. 
\end{cor}

\begin{proof}
The first part is clear from Lemma \ref{f-lem4.1} 
and Lemma \ref{f-lem4.2}. 
This implies that we have maps 
$$
H^q(X, \mathcal C)\to H^q(X, M)\to H^q(X, \mathcal M^{-1}(-S)). 
$$ 
As we saw above, the composition is surjective. 
Hence so is the map on the right. 
\end{proof}

Therefore, 
$H^q(X, \mathcal M^{-1}(-S-D))\to H^q(X, 
\mathcal M^{-1}(-S))$ is 
surjective for every $q$. By Serre duality, we obtain that 
$$H^q(X, \mathcal O_X(K_X)\otimes \mathcal M(S))\to H^q(X, \mathcal 
O_X(K_X)\otimes \mathcal M(S+D))$$ is 
injective for every $q$. 
This means that $$H^q(X, \mathcal O_X(L))\to 
H^q(X, \mathcal O_X(L+D))$$ is injective for every $q$.
\end{proof}

Let us prove Theorem \ref{f-thm1.1}, the main theorem of this paper. 
The proof of Theorem \ref{f-thm1.4} works for Theorem \ref{f-thm1.1} 
with some minor modifications. 

\begin{proof}
[Proof of Theorem \ref{f-thm1.1}] Without loss of generality, we may assume that 
$X$ is connected. 
We can take an effective Cartier divisor $D'$ on $X$ such that 
$D'-D$ is effective and $\Supp D'\subset \Supp \Delta$. 
Therefore, by replacing $D$ with $D'$, we may assume that 
$D$ is a Cartier divisor. 
We set $S=\lfloor \Delta\rfloor$ and $B=\{\Delta\}$. 
By Lemma \ref{f-lem2.11}, we may assume that $B$ is a $\mathbb Q$-divisor.  
We set $\mathcal M=\mathcal O_X(L-K_X-S)$. 
Let $N$ be the smallest positive integer 
such that $N L\sim N(K_X+ S+ B)$. 
We define an $\mathcal O_X$-algebra structure of $\bigoplus _{i=0}^{N-1}
\mathcal M^{-i}(\lfloor iB\rfloor)$ by $s\in H^0(X, \mathcal M^N)$ with 
$(s=0)=NB$. 
We set 
$$\pi:Y=\Spec _X\bigoplus _{i=0}^{N-1}\mathcal M^{-i}(\lfloor iB\rfloor)\to X$$ and 
$T=\pi^*S$. 
Let $Y=\sum _{j\in J}Y_j$ be the irreducible decomposition. 
Then every irreducible component of $Y_{j_1}\cap\cdots \cap Y_{j_l}$ has only 
quotient singularities for every $\{j_1, \cdots, j_l\}\subset J$. 
Let $T=\sum_{i\in I}T_i$ be the irreducible decomposition. 
Then every irreducible component of $T_{i_1}\cap\cdots \cap T_{i_k}$ has only 
quotient singularities for every $\{i_1, \cdots, i_k\}\subset I$. 
Hence it is easy to see that both $Y$ and $T$ are Du Bois 
by Theorem \ref{f-thm3.6} and Lemma \ref{f-lem3.7} (see also \cite{ishida}). 
Therefore, the pair $(Y, T)$ is a Du Bois pair by Proposition \ref{f-prop3.5}. 
This means that 
$\mathcal O_Y(-T)\to 
\underline{\Omega}^0_{Y, T}
$ is a quasi-isomorphism. 
See also \cite[3.4]{ffs}. 
We note that $T$ is Cartier. Hence $\mathcal O_Y(-T)$ is the defining 
ideal sheaf of $T$ on $Y$. 
The $E_1$-degeneration of  
$$
E^{p, q}_1=
\mathbb H^q(Y, \underline {\Omega}^p_{Y, T})\Rightarrow H^{p+q}(Y, j_!\mathbb C_{Y-T})
$$ 
implies that the homomorphism 
$$
H^q(Y, j_!\mathbb C_{Y-T})\to H^q(Y, \mathcal O_Y(-T))
$$
induced by the natural inclusion 
$$
j_!\mathbb C_{Y-T}\subset \mathcal O_Y(-T)
$$ 
is surjective for every $q$ (see Remark \ref{f-rem3.3}). 
By taking a suitable direct summand $$\mathcal C\subset \mathcal M^{-1}(-S)$$ of 
$$\pi_*(j_!\mathbb C_{Y-T})\subset \pi_*\mathcal O_Y(-T), $$ we 
obtain a surjection 
$$
H^q(X, \mathcal C)\to H^q(X, \mathcal M^{-1}(-S)) 
$$ 
induced by the natural inclusion $\mathcal C\subset \mathcal M^{-1}(-S)$ for every $q$. 
It is easy to see that 
Lemma \ref{f-lem4.1} holds for this new setting. Hence Corollary \ref{f-cor4.3} 
also holds without any modifications. 
Therefore, $$H^q(X, \mathcal M^{-1}(-S-D))\to 
H^q(X, \mathcal M^{-1}(-S))$$ is surjective 
for every $q$. By Serre duality, 
we obtain that 
$$
H^q(X, \mathcal O_X(L))\to H^q(X, \mathcal O_X(L+D))
$$ is injective for every $q$. 
\end{proof}

\section{Miscellaneous comments}\label{f-sec5} 

In this section, we collect some miscellaneous comments on related topics. 

\subsection{Ambro's injectivity theorems}\label{f-subsec5.1} 

Let $X$ be a smooth variety and let $\Sigma$ be a simple normal crossing divisor on $X$. 
In order to prove the main theorem of \cite{ambro} (see Theorem \ref{f-thm1.4}), 
Ambro used the complex $(\Omega^\bullet_X(*\Sigma), F)$ and 
the natural inclusion 
$$
(\Omega^\bullet_X(\log \Sigma), F)\subset (\Omega^\bullet_X(*\Sigma), F). 
$$
Hence the arguments in \cite{ambro} 
are different from the proof of Theorem \ref{f-thm1.4} given 
in Section \ref{f-sec4}. 
We do not know how to generalize his approach to the case 
when $X$ is a simple normal crossing variety and $\Sigma$ is a simple normal crossing 
divisor on $X$. 

\subsection{Extension theorem from log canonical centers}\label{f-subsec5.2}
The following result is a slight generalization of \cite[Theorem 6.4]{ambro}. 
Note that \cite[Proposition 5.12]{fujino-gongyo}, which is 
closely related to the abundance conjecture, is a special case of 
Theorem \ref{f-thm5.2.1}. 

\begin{thm-sub}[Extension theorem]\label{f-thm5.2.1}
Let $(X, \Delta)$ be a proper 
log canonical pair. Let $L$ be a Cartier divisor 
on $X$ such that $H=L-(K_X+\Delta)$ is a semi-ample $\mathbb R$-divisor on $X$. 
Let $D$ be an effective $\mathbb R$-divisor on $X$ such that  
$D\sim _{\mathbb R}tH$ for some positive real number $t$ and 
let $Z$ be the union of the log canonical centers of $(X, \Delta)$ contained in $\Supp D$. 
Then the natural restriction map 
$$
H^0(X, \mathcal O_X(L))\to H^0(Z, \mathcal O_Z(L))
$$ 
is surjective. 
\end{thm-sub}

\begin{proof}
Let $f:Y\to X$ be a birational morphism 
from a smooth projective variety $Y$ such that 
$\Exc(f)\cup \Supp f^{-1}_*\Delta$ is a simple normal crossing 
divisor on $Y$. Then we can 
write 
$$
K_Y+\Delta_Y=f^*(K_X+\Delta)+E
$$ 
where $E$ is an effective $f$-exceptional Cartier divisor and 
$\Delta_Y$ is a boundary $\mathbb R$-divisor. 
Without loss of generality, we may further 
assume that 
$f^{-1}(Z)$ is a divisor on $Y$. Let $W$ be the union of all the 
log canonical centers of $(Y, \Delta_Y)$ whose images by $f$ 
are contained in $Z$. Note that $W$ is a divisor on $Y$ such that 
$W\leq \lfloor \Delta_Y\rfloor$. 
We consider the short exact sequence 
$$
0\to \mathcal O_Y(E-W)\to \mathcal O_Y(E)\to \mathcal O_W(E)\to 0. 
$$ 
Since 
$$
E-W=K_Y+(\Delta_Y-W)-f^*(K_X+\Delta), 
$$ 
there are no associated primes of $R^1f_*\mathcal O_Y(E-W)$ in 
$Z=f(W)$ by \cite[Theorem 6.3 (i)]{fujino-fund}. 
Therefore, the connecting homomorphism 
$$
\delta: f_*\mathcal O_W(E)\to R^1f_*\mathcal O_Y(E-W)
$$ 
is zero. Hence we obtain 
$$
\mathcal O_X\simeq f_*\mathcal O_Y(E)\to f_*\mathcal O_W(E)
$$ 
is surjective. This implies that $f_*\mathcal O_W(E)\simeq \mathcal O_Z$. 
Since $H^0(Y, \mathcal O_Y(f^*L+E))\simeq H^0(X, \mathcal O_X(L))$ and 
$H^0(W, \mathcal O_W(f^*L+E))\simeq H^0(Z, \mathcal O_Z(L))$, it is sufficient 
to prove that the natural restriction map 
$$
H^0(Y, \mathcal O_Y(f^*L+E))\to H^0(W, \mathcal O_W(f^*L+E))
$$ 
is surjective. By assumption, there is a morphism $g:X\to V$ such that 
$V$ is a normal projective variety, 
$g_*\mathcal O_X\simeq \mathcal O_V$, 
and $H\sim _{\mathbb R}g^*A$, where 
$A$ is an ample $\mathbb R$-divisor on $V$. 
We note that 
\begin{align*}
(f^*L+E-W)-(K_Y+\Delta_Y-W)&=f^*(L-(K_X+\Delta))\\
&\sim _{\mathbb R}f^*g^*A. 
\end{align*} 
By the assumption on $D$ and the construction of 
$$
Y\overset{f}\longrightarrow X\overset{g}\longrightarrow V, 
$$ 
we can find an effective ample Cartier divisor $D_1$ and an effective 
ample $\mathbb R$-divisor $D_2$ on $V$ such that 
$D_1+D_2\sim _{\mathbb R}sA$ for some positive real number $s$, 
$W\leq f^*g^*D_1$, and that $\Supp f^*g^*(D_1+D_2)$ 
contains no log canonical centers of $(Y, \Delta_Y-W)$. Hence 
$$
H^i(Y, \mathcal O_Y(f^*L+E-W))\to H^i(Y, \mathcal O_Y(f^*L+E)) 
$$ 
is injective for every $i$ (see \cite[Theorem 6.1]{fujino-fund}). 
See also Theorem \ref{f-thm1.5}.  
In particular, 
$$
H^1(Y, \mathcal O_Y(f^*L+E-W))\to H^1(Y, \mathcal O_Y(f^*L+E)) 
$$ 
is injective. Thus we obtain that 
$$
H^0(Y, \mathcal O_Y(f^*L+E))\to H^0(W, \mathcal O_W(f^*L+E)) 
$$ is surjective. 
Therefore, we obtain the desired surjection. 
\end{proof}

The proof of Theorem \ref{f-thm5.2.1} is essentially 
the same as that of \cite[Proposition 5.12]{fujino-gongyo} 
and is different from the arguments in \cite[Section 6]{ambro}. 
The framework discussed in \cite{fujino-fund} is sufficient for 
Theorem \ref{f-thm5.2.1}. 
We recommend the reader to compare 
the above proof with the proof of \cite[Theorem 6.4]{ambro}, 
which is shorter than our proof and 
is based on \cite[Theorem 6.2]{ambro}. 
We will present the original proof of Theorem \ref{f-thm5.2.1} 
as an application of Theorem \ref{f-thm5.3.3} below for the 
reader's convenience. For the relative version of Theorem \ref{f-thm5.2.1}, 
see Theorem \ref{f-thm6.4} below. 

\subsection{The maximal non-lc ideal sheaves}\label{f-subsec5.3}
By combining Theorem \ref{f-thm1.4} with the notion of {\em{maximal}} 
non-lc ideal sheaves, 
we have some interesting results due to Ambro (\cite{ambro}). 
Note that the ideal sheaf defined in \cite[Definition 4.3]{ambro} 
is nothing but the {\em{maximal non-lc ideal sheaf}} introduced 
in \cite[Definition 7.1]{fst} (see also \cite[Remark 7.6]{fujino-fund}). 

Let us recall the definition of maximal non-lc ideal sheaves. 

\begin{defn-sub}[Maximal non-lc ideal sheaves]\label{f-def5.3.1} 
Let $X$ be a normal variety and let $\Delta$ be an 
$\mathbb R$-divisor on $X$ such that $K_X+\Delta$ is $\mathbb R$-Cartier. 
Let $f:Y\to X$ be a resolution 
with $$K_Y+\Delta_Y=f^*(K_X+\Delta)$$ such that 
$\Supp \Delta_Y$ is a simple normal crossing divisor. 
Then we put 
$$
\mathcal J'(X, \Delta)=f_*\mathcal O_Y(\lceil K_Y-f^*(K_X+\Delta)
+\varepsilon F\rceil)
$$ 
for $0<\varepsilon \ll 1$, where $F=\Supp \Delta^{\geq 1}_Y$. 
We call $\mathcal J'(X, \Delta)$ the {\em{maximal non-lc 
ideal sheaf}} associated to 
$(X, \Delta)$. 
It is easy to see that 
$$
\mathcal J'(X, \Delta)=f_*\mathcal O_Y(-\lfloor \Delta_Y\rfloor+\sum _{k=1}^\infty 
{}^k\!\Delta_Y).  
$$ 
Note that there is a positive integer $k_0$ such that ${}^k\!\Delta_Y=0$ for every 
$k>k_0$. Therefore, 
$$\sum _{k=1}^\infty{}^k\!\Delta_Y={}^1\!\Delta_Y+
{}^2\!\Delta_Y+\cdots +{}^{k_0}\!\Delta_Y. $$  
We also note that 
$$
\mathcal J_{\mathrm{NLC}}(X, \Delta)=f_*\mathcal O_Y(-\lfloor \Delta_Y\rfloor+
\Delta^{=1}_Y)
$$ 
is the 
({\em{minimal}}) {\em{non-lc ideal sheaf}} associated to $(X, \Delta)$ 
and  that 
$$
\mathcal J(X, \Delta)=f_*\mathcal O_Y(-\lfloor \Delta_Y\rfloor)
$$ 
is the {\em{multiplier ideal sheaf}} associated to $(X, \Delta)$. 
It is obvious that 
$$
\mathcal J(X, \Delta)\subset \mathcal J_{\mathrm{NLC}}(X, \Delta) 
\subset\mathcal J'(X, \Delta). 
$$
\end{defn-sub}

For the 
details of $\mathcal J'(X, \Delta)$, see \cite{fst} (see also \cite{fujino-nonlc}). 

\begin{rem-sub}[Non-F-pure ideals]\label{f-rem5.3.2}
A positive characteristic analog of 
$\mathcal J'(X, \Delta)$, which we call the {\em{non-F-pure ideal}} 
associated to $(X, \Delta)$ 
and is denoted by $\sigma(X, \Delta)$, introduced 
in \cite{fst} is now becoming a very important tool 
for higher-dimensional algebraic geometry in positive characteristic. 
\end{rem-sub}

Theorem \ref{f-thm5.3.3} is a nontrivial application of 
Theorem \ref{f-thm1.4}. For the relative version of 
Theorem \ref{f-thm5.3.3}, see Theorem \ref{f-thm6.2} below. 

\begin{thm-sub}[{\cite[Theorem 6.2]{ambro}}]\label{f-thm5.3.3}  
Let $X$ be a proper normal variety and let $\Delta$ be 
an effective $\mathbb R$-divisor on $X$ such that 
$K_X+\Delta$ is $\mathbb R$-Cartier. 
Let $L$ be a Cartier divisor on $X$ such that 
$L-(K_X+\Delta)$ is semi-ample. 
Let $\mathcal J'(X, \Delta)$ be the maximal 
non-lc ideal sheaf associated to $(X, \Delta)$ and let $Y$ be the closed 
subscheme defined by $\mathcal J'(X, \Delta)$. Then we 
have a short exact sequence 
\begin{align*}
0&\to H^0(X, \mathcal J'(X, \Delta)\otimes \mathcal O_X(L))
\\ &\to H^0(X, \mathcal O_X(L))
\to H^0(Y, \mathcal O_Y(L))\to 0. 
\end{align*}
\end{thm-sub}

We describe the proof of Theorem \ref{f-thm5.3.3} for the reader's convenience 
(see also \cite{ambro}). 

\begin{proof}
We take an effective general $\mathbb R$-divisor $D$ with small coefficients such that 
$L-(K_X+\Delta)\sim _{\mathbb R}D$. By replacing $\Delta$ with 
$\Delta+D$, we may assume that 
$L\sim _{\mathbb R}K_X+\Delta$. 
Let $Z\to X$ be a resolution such that 
$K_Z+\Delta_Z=f^*(K_X+\Delta)$. 
We may assume that $\Supp \Delta_Z$ is 
a simple normal crossing divisor. 
We note that 
$$
-\lfloor \Delta_Z\rfloor+\sum _{k=1}^\infty 
{}^k\!\Delta_Z=(K_Z+\{\Delta_Z\}+\sum _{k=1}^\infty{}^k\!\Delta_Z)-f^*(K_X+\Delta). 
$$ 
We write 
$$
-\lfloor \Delta_Z\rfloor+\sum _{k=1}^\infty 
{}^k\!\Delta_Z=P-N
$$ 
where $P$ and $N$ are effective and have no common irreducible components. 
Note that $P$ is $f$-exceptional since $\Delta$ is effective. 
Therefore, 
$$
f^*L+P-N\sim _{\mathbb R}K_Z+\{\Delta_Z\}+\sum _{k=1}^\infty 
{}^k\!\Delta_Z. 
$$ 
Thus 
$$
H^i(Z, \mathcal O_Z(f^*L+P-N))\to H^i(Z, \mathcal O_Z(f^*L+P))
$$ 
is injective for every $i$ by Theorem \ref{f-thm1.4}. 
This is because $$\Supp N\subset \Supp (\{\Delta_Z\}+\sum _{k=1}^\infty 
{}^k\!\Delta_Z).$$ 
We note that $$f_*\mathcal O_Z(f^*L+P-N)\simeq \mathcal J'(X, \Delta)\otimes 
\mathcal O_X(L)$$ and $$f_*\mathcal O_Z(f^*L+P)\simeq \mathcal O_X(L).$$ 
By the following commutative diagram: 
$$
\xymatrix{
H^1(Z, \mathcal O_Z(f^*L+P-N))
\ar[r]^{\quad b}& 
H^1(Z, \mathcal O_Z(f^*L+P)) \\
H^1(X, \mathcal J'(X, \Delta)\otimes \mathcal O_X(L))\ar[r]_{\quad \quad d}
\ar[u]^{a}&
H^1(X, \mathcal O_X(L))
\ar[u]_{c}, 
}
$$ 
we obtain that $$H^1(X, \mathcal J'(X, \Delta)\otimes \mathcal O_X(L))
\to H^1(X, \mathcal O_X(L))$$ 
is injective. Note that $a$ and $c$ are injective by the Leray spectral 
sequences and that $b$ is injective by the above argument. 
Hence the natural restriction map 
$$
H^0(X,  \mathcal O_X(L))\to H^0(Y, \mathcal O_Y(L)) 
$$ 
is surjective. 
We obtain the desired short exact sequence. 
\end{proof}

Theorem \ref{f-thm5.3.3} shows 
that $\mathcal J'(X, \Delta)$ is useful for some applications. 
We give the original proof of Theorem \ref{f-thm5.2.1} as an application of 
Theorem \ref{f-thm5.3.3}. 

\begin{proof}[{Proof of Theorem \ref{f-thm5.2.1}}] 
Let $\varepsilon$ be a small positive number. 
Then it is easy to see that 
$\mathcal J'(X, \Delta+\varepsilon D)=\mathcal I_Z$, where 
$\mathcal I_Z$ is the defining ideal sheaf of 
$Z$. Since $L-(K_X+\Delta+\varepsilon D)\sim _
{\mathbb R}(1-\varepsilon t)H$ is semi-ample, 
we have the following short exact sequence 
\begin{align*}
0\to H^0(X, \mathcal J'(X, \Delta+\varepsilon D)\otimes \mathcal O_X(L))
&\to H^0(X, \mathcal O_X(L))\\
&\to H^0(Z, \mathcal O_Z(L))\to 0
\end{align*}
by Theorem \ref{f-thm5.3.3}. In particular, the natural restriction map 
$$
H^0(X, \mathcal O_X(L))\to H^0(Z, \mathcal O_Z(L))
$$ 
is surjective. 
\end{proof}

The following theorem is Ambro's inversion of adjunction. 
For the relative version of Theorem \ref{f-thm5.3.4}, see Theorem \ref{f-thm6.3} below. 

\begin{thm-sub}[{\cite[Theorem 6.3]{ambro}}]\label{f-thm5.3.4}
Let $X$ be a proper normal irreducible variety and let 
$\Delta$ be an effective $\mathbb R$-divisor on $X$ such that 
$-(K_X+\Delta)$ is semi-ample. 
Suppose that the non-lc locus $\Nlc(X, \Delta)$ of $(X, \Delta)$ is not empty, that is, 
$(X, \Delta)$ is not log canonical. 
Then $\Nlc (X, \Delta)$ is connected and intersects every 
log canonical center of $(X, \Delta)$. 
\end{thm-sub}

We describe Ambro's proof of Theorem \ref{f-thm5.3.4} based 
on Theorem \ref{f-thm1.4} in order to show how to use Theorem \ref{f-thm1.4}. 

\begin{proof} 
We take an effective general $\mathbb R$-divisor $D$ with small coefficients such that 
$D\sim _{\mathbb R}-(K_X+\Delta)$. 
By replacing $\Delta$ with $\Delta+D$, we may assume that 
$K_X+\Delta\sim _{\mathbb R}0$. 
We set $Y=\Nlc (X, \Delta)$. 
By Theorem \ref{f-thm5.3.3}, we 
have the following short exact sequence:  
$$
0\to H^0(X, \mathcal J'(X, \Delta))\to H^0(X, \mathcal O_X)\to 
H^0(Y, \mathcal O_Y)\to 0. 
$$ 
This implies that $H^0(Y, \mathcal O_Y)\simeq \mathbb C$. 
Hence $Y$ is connected. Let $C$ be a log canonical 
center of $(X, \Delta)$. 
Let $f:Z\to X$ be a resolution such that $\Exc(f)\cup \Supp f^{-1}_*\Delta$ is 
a simple normal crossing divisor and that 
$f^{-1}(C)$ is a divisor. 
We set $K_Z+\Delta_Z=f^*(K_X+\Delta)$. 
Let $W$ be the union of all 
the irreducible components of $\Delta^{=1}_Z$ whose 
images by $f$ are contained in $C$. 
It is obvious that $f(W)=C$. 
By construction, 
we have 
$$
-\lfloor \Delta_Z\rfloor +\sum _{k=1}^{\infty}{}^k\!\Delta_Z-W\sim _{\mathbb R}
K_Z+\{\Delta_Z\}+\sum _{k=1}^{\infty}{}^k\!\Delta_Z-W
$$ 
since $K_Z+\Delta_Z\sim _{\mathbb R}0$. 
We set 
$$
-\lfloor \Delta_Z\rfloor +\sum _{k=1}^{\infty}{}^k\!\Delta_Z=P-N
$$ 
where $P$ and $N$ are effective and have no common irreducible components. 
Note that $P$ is $f$-exceptional. 
By Theorem \ref{f-thm1.4}, 
$$
H^i(Z, \mathcal O_Z(P-N-W))\to H^i(Z, \mathcal O_Z(P-W))
$$ 
is injective for every $i$ because $\Supp N\subset \Supp 
(\{\Delta_Z\}+\sum _{k=1}^\infty{}^k\!\Delta_Z-W)$. 
Thus the natural restriction map 
$$
H^0(Z, \mathcal O_Z(P-W))\to H^0(N, \mathcal O_N(P-W))
$$ 
is surjective. Since $H^0(Z, \mathcal O_Z(P-W))=0$, 
we obtain $H^0(N, \mathcal O_N(P-W))=0$. 
On the other hand, 
$$H^0(N, \mathcal O_N(P-W))\subset H^0(N, \mathcal O_N(P))\ne 0
$$ 
implies $N\cap W\ne \emptyset$. Thus we obtain 
$C\cap Y\ne \emptyset$. 
\end{proof}

\begin{rem-sub}\label{f-rem5.3.5}
If $X$ is {\em{projective}} in Theorem \ref{f-thm5.3.4}, then we can 
prove Theorem \ref{f-thm5.3.4} without using Theorem \ref{f-thm5.3.3}. 
We give a sketch of the proof. 
We may assume that $K_X+\Delta\sim _{\mathbb R}0$. 
Let $f:Y\to X$ be a dlt blow-up 
with $K_Y+\Delta_Y=f^*(K_X+\Delta)$. 
We may assume that 
$a(E, X, \Delta)\leq -1$ for every $f$-exceptional 
divisor  and that 
$(Y, \Delta^{\leq 1}_Y+S)$ is a dlt pair where $S=\Supp \Delta^{>1}_Y$. We 
run a minimal model program with respect to 
$K_Y+\Delta^{\leq 1}_Y+S$. 
Note that $K_Y+\Delta^{\leq 1}_Y+S\sim _{\mathbb R} S-\Delta^{>1}_Y\neq 0$ 
is not pseudo-effective.
By the similar argument to the proof of \cite[Proposition 2.1]
{fujino-abundance} (cf.~\cite[Theorem 3.47]{book}), 
we can recover Theorem \ref{f-thm5.3.4} when $X$ is projective. We leave the details 
as exercises for the interested reader. 
\end{rem-sub}

\section{Relative version}\label{f-sec6}
In this section, we discuss the relative version of Theorem \ref{f-thm1.1} 
and some related results. 

\begin{thm}[Relative injectivity theorem]\label{f-thm6.1}
Let $X$ be a simple normal crossing variety and let $\Delta$ be an 
$\mathbb R$-Cartier $\mathbb R$-divisor 
on $X$ such that $\Supp \Delta$ is a simple normal crossing divisor on $X$ 
and that $\Delta$ is a boundary $\mathbb R$-divisor on $X$. 
Let $\pi:X\to V$ be a proper morphism between algebraic varieties 
and let $L$ be a Cartier divisor 
on $X$ and let $D$ be an effective Weil divisor on $X$ whose support 
is contained in $\Supp \Delta$. 
Assume that 
$L\sim _{\mathbb R, \pi}K_X+\Delta$, that is, 
there is an $\mathbb R$-Cartier divisor $B$ on $V$ such that 
$L\sim _\mathbb R K_X+\Delta+\pi^*B$. 
Then the natural homomorphism 
$$
R^q\pi_*\mathcal O_X(L)\to R^q\pi_*\mathcal O_X(L+D)
$$ 
induced by the inclusion $\mathcal O_X\to \mathcal O_X(D)$ is injective for every $q$. 
\end{thm}

By using \cite{bierstone-p} (see 
\cite[Lemma 3.6]{fujino-vani}), 
we can reduce Theorem \ref{f-thm6.1} to Theorem \ref{f-thm1.1}. 

\begin{proof}
By shrinking $V$, we may assume that $V$ is affine and $L\sim _{\mathbb R}K_X+\Delta$. 
Without loss of generality, we may assume that $X$ is connected. 
Let $\overline V$ be a projective compactification of $V$. 
By \cite[Lemma 3.6]{fujino-vani}, we 
can compactify $\pi:X\to V$ to $\overline {\pi}: \overline X\to \overline V$. 
By the same argument as in Step 2 in the proof of \cite[Theorem 3.7 (i)]{fujino-vani}, 
we may assume that there is a Cartier divisor $\overline L$ on $\overline X$ such that 
$\overline L|_X=L$. We can write 
$$
L-(K_X+\Delta)=\sum _i b_i (f_i)
$$
where $b_i$ is a real number and $f_i\in \Gamma (X, \mathcal K^*_X)$ for 
every $i$. 
We put 
$$
E=\sum _i b_i (f_i)-(\overline L-(K_{\overline X}+\overline \Delta)). 
$$ 
Then we have $$
\overline L+\lceil E\rceil \sim _{\mathbb R}K_{\overline X}+\overline \Delta+\{-E\}. 
$$
By the above construction, it is obvious that $\Supp E\subset \overline X\setminus X$. 
Let $\overline D$ be the closure of $D$ in $\overline X$. 
It is sufficient to prove that the map
$$
\varphi^q: R^q\overline {\pi}_*\mathcal O_{\overline X}(\overline {L}+\lceil E\rceil)\to 
R^q\overline{\pi}_*\mathcal O_{\overline X}(\overline {L}+\lceil E\rceil+\overline D)
$$ 
induced by the natural inclusion $\mathcal O_{\overline X}\to \mathcal O_{\overline X}
(\overline D)$ is injective for every $q$.  Suppose that $\varphi^q$ is not 
injective for some $q$. Let $A$ be a 
sufficiently ample general Cartier divisor on $\overline V$ 
such that $H^0(\overline V, 
\mathrm{Ker}{\varphi^q}\otimes \mathcal O_{\overline V}(A))\ne 0$. 
In this case, the map
\begin{align*}
&H^0(\overline V, R^q\overline {\pi}_*\mathcal O_{\overline X}(\overline {L}+\lceil E\rceil)
\otimes \mathcal O_{\overline V}(A))\\ &\to 
H^0(\overline V, R^q\overline {\pi}_*\mathcal O_{\overline X}(\overline {L}+\lceil E\rceil
+\overline D)
\otimes \mathcal O_{\overline V}(A))
\end{align*} 
induced by $\varphi^q$ 
is not injective. Since $A$ is sufficiently ample, this implies 
that 
\begin{align*}
&H^q(\overline X, \mathcal O_{\overline X}(\overline {L}+\lceil E\rceil+\overline {\pi}^*A))
\\&\to 
H^q(\overline X, 
\mathcal O_{\overline X}(\overline {L}+\lceil E\rceil+\overline {\pi}^*A+\overline D))
\end{align*} 
is not injective. Since 
$$
\overline L+\lceil E\rceil +\overline{\pi}^*A 
\sim _{\mathbb R}K_{\overline X}+\overline \Delta+\{-E\} 
+\overline{\pi}^*A,  
$$ 
it contradicts Theorem \ref{f-thm1.1}. 
Hence $\varphi^q$ is injective for every $q$. 
\end{proof}

The following theorem is the relative version of Theorem \ref{f-thm5.3.3}. 
It is obvious by the proof of Theorem \ref{f-thm5.3.3} and 
Theorem \ref{f-thm6.1}. 

\begin{thm}\label{f-thm6.2}
Let $X$ be a normal variety and let $\Delta$ be 
an effective $\mathbb R$-divisor on $X$ such that 
$K_X+\Delta$ is $\mathbb R$-Cartier. 
Let $\pi:X\to V$ be a proper morphism between algebraic varieties 
and let $L$ be a Cartier divisor on $X$ such that 
$L-(K_X+\Delta)$ is semi-ample over $V$. 
Let $\mathcal J'(X, \Delta)$ be the maximal 
non-lc ideal sheaf associated to $(X, \Delta)$ and let $Y$ be the closed 
subscheme defined by $\mathcal J'(X, \Delta)$. Then we 
have a short exact sequence 
\begin{align*}
0\to \pi_*(\mathcal J'(X, \Delta)\otimes \mathcal O_X(L))\to \pi_*\mathcal O_X(L)
\to \pi_*\mathcal O_Y(L)\to 0. 
\end{align*}
\end{thm}
\begin{proof}
It is sufficient to 
prove that $\pi_*\mathcal O_X(L)\to \pi_*\mathcal O_Y(L)$ is surjective. 
Since the problem is local, we may assume that $V$ is affine by shrinking 
$V$. Then the proof of Theorem \ref{f-thm5.3.3} works 
without any modifications if we use Theorem \ref{f-thm6.1}. 
\end{proof}

The relative version of Theorem \ref{f-thm5.3.4} is: 

\begin{thm}\label{f-thm6.3}
Let $X$ be a normal variety and let $\pi:X\to V$ be 
a proper morphism between algebraic 
varieties with $\pi_*\mathcal O_X\simeq \mathcal O_V$. 
Let 
$\Delta$ be an effective $\mathbb R$-divisor on $X$ such that 
$-(K_X+\Delta)$ is semi-ample over $V$. Let $x$ be a closed point of $V$. 
Suppose that $$\Nlc(X, \Delta)\cap\pi^{-1}(x)\ne \emptyset. $$ 
Then $\Nlc (X, \Delta)\cap \pi^{-1}(x)$ is connected and intersects every 
log canonical center $C$ of $(X, \Delta)$ with $C\cap \pi^{-1}(x)\ne \emptyset$. 
\end{thm}

\begin{proof}
By shrinking $V$, we may assume that $V$ is affine. 
As in the proof of Theorem \ref{f-thm5.3.4}, we may assume that 
$K_X+\Delta\sim _{\mathbb R}0$. 
From now on, we use the same notation as in the proof of 
Theorem \ref{f-thm5.3.4}. 
Since $\mathcal O_V\simeq \pi_*\mathcal O_X\to \pi_*\mathcal O_Y$ is surjective 
by Theorem \ref{f-thm6.2}, 
$Y\cap \pi^{-1}(x)$ is connected. 
By Theorem \ref{f-thm6.1}, 
$$
R^i(\pi\circ f)_*\mathcal O_Z(P-N-W)\to 
R^i(\pi\circ f)_*\mathcal O_Z(P-W)
$$ is injective for every $i$. 
Thus the natural restriction map  
$$
(\pi\circ f)_*\mathcal O_Z(P-W)\to (\pi\circ f)_*\mathcal O_N(P-W)
$$ is surjective. 
Since $(\pi\circ f)_*\mathcal O_Z(P-W)\subset \mathcal I_x
\subsetneq \mathcal O_V$, where $\mathcal I_x$ is the defining ideal sheaf of $x$ 
on $V$, 
we obtain 
$$(\pi\circ f)_*\mathcal O_N(P-W)\subsetneq (\pi\circ f)_*\mathcal O_N
\subset (\pi\circ f)_*\mathcal O_N(P)$$ at $x$.  
This implies $N\cap W\cap (\pi\circ f)^{-1}(x)\ne \emptyset$. 
Therefore, $C\cap Y\cap \pi^{-1}(x)\ne \emptyset$. 
\end{proof}

Theorem \ref{f-thm6.4}, which is the relative version of Theorem \ref{f-thm5.2.1}, 
directly follows from Theorem \ref{f-thm6.2}. 
See the proof of Theorem \ref{f-thm5.2.1} by 
Theorem \ref{f-thm5.3.3} in Subsection \ref{f-subsec5.3}. 

\begin{thm}[Relative extension theorem]\label{f-thm6.4}
Let $(X, \Delta)$ be a 
log canonical pair and let $\pi:X\to V$ be 
a proper morphism. Let $L$ be a Cartier divisor 
on $X$ such that $H=L-(K_X+\Delta)$ is a $\pi$-semi-ample $\mathbb R$-divisor on $X$. 
Let $D$ be an effective $\mathbb R$-divisor on $X$ such that  
$D\sim _{\mathbb R, \pi}tH$, that is, 
there is an $\mathbb R$-Cartier divisor $B$ on $V$ with 
$D\sim _\mathbb R tH+\pi^*B$,  
for some positive real number $t$ and 
let $Z$ be the union of the log canonical centers of $(X, \Delta)$ contained in $\Supp D$. 
Then the natural restriction map 
$$
\pi_*\mathcal O_X(L)\to \pi_*\mathcal O_Z(L)
$$ 
is surjective. 
\end{thm}


\end{document}